\documentclass[a4paper, BCOR=0.0mm, DIV=calc]{scrartcl}
\usepackage[T1]{fontenc}
\usepackage[utf8]{inputenc}
\usepackage[ngerman]{babel}
\usepackage[osf,sc]{mathpazo}
\usepackage{textcomp}
\usepackage{microtype}
\usepackage{amsmath, amssymb, amsfonts}
\usepackage{tabularx}
\usepackage{nicefrac}
\KOMAoptions{twoside=false, twocolumn=false, headinclude=false, footinclude=false, mpinclude=false, pagesize=auto}
\recalctypearea
\newcommand{\ebinom}[2]{\left(\frac{#1}{#2} \right)}
\setkomafont{section}{\mdseries\scshape\Large}
\setkomafont{title}{\mdseries\scshape}

\begin{document}
\title{Ein Kommentar zum Kettenbruch, mit welchem der bedeutende Lagrange die Binomialpotenzen ausgedrückt hat\footnote{
Originaltitel: "`Commentatio in fractionem continuam, qua illustris La Grange potestates binomiales expressit"', erstmals publiziert in "`\textit{Mémoires de l'académie des sciences de St.-Petersbourg} 6, 1818, pp. 3-11"', Nachdruck in "`\textit{Opera Omnia}: Series 1, Volume 16, pp. 232 - 240"', Eneström-Nummer E750, übersetzt von: Alexander Aycock, Textsatz: Artur Diener,  im Rahmen des Projektes "`Eulerkreis Mainz"' }}
\author{Leonhard Euler}
\date{}
\maketitle
\paragraph{§1}
Dieser bedeutende Mann hat diese Binomialpotenz $(1+x)^n$ durch eine völlig einzigartige Methode aus seinem logarithmischen Differential in diesen Kettenbruch verwandelt:
\[
	(1+x)^n = 1 + \cfrac{nx}{1 + \cfrac{(1-n)x}{2 + \cfrac{(1+n)x}{3 + \cfrac{(2-n)x}{2 + \cfrac{(2+n)x}{5 + \cfrac{(3-n)x}{2 + \cfrac{(3+n)x}{7 + \mathrm{etc}}}}}}}}
\]
welcher Ausdruck sich der hervorstechenden Eigenschaft erfreut, dass er, sooft der Exponent $n$ eine ganze Zahl war, ob positiv oder negativ, abbricht und auf endliche Form zurückgeht.
\paragraph{§2}
Weil ja dieser Kettenbruch nicht nach einem gleichmäßigen Gesetz, sondern nach einem unterbrochenen, voranschreitet, wollen wir ihn auf ein gleichmäßiges Bildungsgesetz zurückführen; das kann am angenehmsten gemacht werden, wenn wir ihn auf die folgende Weise durch Teile darstellen:
\begin{align*}
	(1+x)^n &= 1 + \frac{nx}{A} \\[1mm]
	A &= 1 + \cfrac{(1-n)x}{2 + \cfrac{(1+n)x}{B}} \\[1mm]
	B &= 1 + \cfrac{(2-n)x}{2 + \cfrac{(2+n)x}{C}} \\[1mm]
	C &= 1 + \cfrac{(3-n)x}{2 + \cfrac{(3+n)x}{D}} \\[1mm]
	D &= 1 + \cfrac{(4-n)x}{2 + \cfrac{(4+n)x}{E}} \\[1mm]
	\mathrm{etc.}
\end{align*}
Daher werden wir also durch Reduktion haben:
\begin{align*}
A = 1 + \frac{(1-n)Bx}{2B + (1+n)x} &= 1 + \frac{(1-n)x}{2} - \frac{(1-nn)xx : 2}{2B + (1+n)x} \\
	&= 1 + \frac{(1-n)x}{2} + \frac{(nn-1)xx : 4}{B + \ebinom{1+n}{2}x}.
\end{align*}
Auf ähnliche Weise wird
\begin{align*}
B = 3 + \frac{(2-n)Cx}{2C + (2+n)x} &= 3 + \frac{(2-n)x}{2} - \frac{(4-nn)xx : 2}{2C + (2+n)x} \\
 &= 3 + \frac{(2-n)x}{2} + \frac{(nn-4)xx : 4}{C + \ebinom{2+n}{2}x}
\end{align*}
sein. Auf dieselbe Weise werden wir
\begin{align*}
	C = 5 + \frac{(3-n)Dx}{2D + (3+n)x} &= 5 + \frac{(3-n)x}{2} - \frac{(9-nn)xx : 2}{2D + (3+n)x} \\
	&=5 + \frac{(3-n)x}{2} + \frac{(nn-9)xx : 4}{D + \ebinom{3+n}{2}x}
\end{align*}
haben und so weiter.
\paragraph{§3}
Wenn wir also diese Werte der Reihe nach anstelle von $A$, $B$, $C$ einsetzen, wird der Kettenbruch die folgende Form annehmen:
\[
	(1+x)^n = 1 + \cfrac{nx}{1 + \tfrac{(1-n)x}{2} + \cfrac{(nn-1)xx : 4}{3(1 + \tfrac{1}{2}x) + \cfrac{(nn-4)xx : 4}{5(1+\tfrac{1}{2}x) + \cfrac{(nn-9)xx:4}{7(1+\tfrac{1}{2}x) + \cfrac{(nn-16)xx:4}{\mathrm{etc}}}}}}
\]
\paragraph{§4}
Damit wir hier die Partialbrüche loswerden, wollen wir $x=2y$ setzen, dass wir diesen Ausdruck erhalten:
\[
	(1+2y)^n = 1 + \cfrac{2ny}{1 + (1-n)y + \cfrac{(nn-1)yy}{3(1+y) + \cfrac{(nn-4)yy}{5(1+y) + \cfrac{(nn-9)yy}{7(1+y) + \mathrm{etc}}}}}
\]
welche leicht in diesen verwandelt wird:
\[
	\frac{2ny}{(1+2y)^n - 1} = 1 + (1-n)y + \cfrac{(nn-1)yy}{3(1+y) + \cfrac{(nn-4)yy}{5(1+y) +\mathrm{etc}}}
\]
Man addiere auf beiden Seiten $ny$, dass er 
\[
	\frac{ny(1+(1+2y)^n)}{(1+2y)^n -1} = 1 + y + \cfrac{(nn-1)yy}{3(1+y) + \cfrac{(nn-4)yy}{5(1+y) + \mathrm{etc}}}
\]
welcher Ausdruck schon nach einer hinreichend regelmäßigen Struktur voranschreitet.
\paragraph{§5}
Wir wollen gleich auf beiden Seiten durch $1+y$ teilen, und die linke Seite wird
\[
	\frac{ny}{1+y} \cdot \frac{(1+2y)^n + 1}{(1+2y)^n - 1}
\]
werden. Auf der rechten Seite aber teile man die einzelnen Brüche oben und unten durch $1+y$ und es wird diese Form hervorgehen
\[
	1 + \cfrac{(nn-1)yy : (1+y)^2}{3 + \cfrac{(nn-4)yy : (1+y)^2}{5 + \cfrac{(nn-9)yy : (1+y)^2}{7 + \cfrac{(nn-16)yy : (1+y)^2}{9 + \cfrac{(nn-25)yy : (1+y)^2}{11 + \mathrm{etc}}}}}}.
\]
\paragraph{§6}
Diesen Ausdruck wollen wir also erneut vereinfachen, indem wir $\frac{y}{1+y} = z$ setzen, sodass $y = \frac{z}{1-z}$ ist. Auf diese Weise aber wird diese linke Seite wegen
\[
	1 + 2y = \frac{1+z}{1-z}
\]
diese Form annehmen:
\[
	\frac{nz[(1+z)^n + (1-z)^n]}{(1+z)^n - (1-z)^n},
\]
was also diesem Kettenbruch gleich werden wird:
\[
	1 + \cfrac{(nn-1)zz}{3 + \cfrac{(nn-4)zz}{5 + \cfrac{(nn-9)zz}{7 + \cfrac{(nn-16)zz}{9 + \mathrm{etc}}}}}
\]
welcher, wegen seiner Eleganz, die höchste Aufmerksamkeit verdient.
\paragraph{§7}
Nun ist also per se klar, dass dieser Ausdruck immer irgendwann abbricht, sooft $n$ eine ganze Zahl war, ob positiv oder negativ. Es ist aber evident, dass auch die linke Seite denselben Wert beibehält, auch wenn $-n$ für $n$ geschrieben wird. Nachdem das nämlich gemacht worden ist, wird
\[
	\frac{-nz[(1+z)^{-n} + (1-z)^{-n}]}{(1+z)^{-n} - (1-z)^{-n}}
\]
werden, welcher Bruch, wenn er mit $(1-zz)^n$ erweitert wird, diese Form annehmen wird:
\[
	\frac{-nz[(1-z)^n + (1+z)^n]}{(1-z)^n - (1+z)^n} = \frac{nz[(1+z)^n + (1-z)^n]}{(1+z)^n - (1-z)^n},
\]
was der vorhergehende Ausdruck selbst ist. Und genauso ist es, ob nun dem Buchstaben $n$ ein positiver oder negativer Wert zugeteilt wird.
\paragraph{§8}
Wenn wir so $n=\pm 1$ nehmen, wird die linke Seite gleich $1$, welcher auch der Wert der rechten ist. Weiter wird für $n=\pm 2$ gesetzt die linke Seite gleich $1+zz$, die rechte Seite wird in der Tat auch gleich $1+zz$. Auf ähnliche Weise wird für $n = \pm 3$ genommen die linke Seite, wie auch die rechte, $\frac{3(1+3zz)}{3+zz}$.
\paragraph{§9}
Daher lassen sich hoffentlich aber einige Schlussfolgerungen von größter Bedeutung ableiten, je nachdem ob dem Exponenten $n$ entweder ein verschwindender oder ein unendlicher Wert zugeteilt wird; besonders aber der Fall, in dem dem Buchstaben $z$ ein imaginärer Wert gegeben wird, führt zu einer außerordentlichen Schlussfolgerung, weil ja dieser Kettenbruch selbst nichtsdestoweniger reell bleibt, von welcher Schlussfolgerung aus wir also anfangen wollen.
\section*{Schlussfolgerung 1, \\ bei der $z = t\sqrt{-1}$ ist}
\paragraph{§10}
In diesem Fall wird also dieser Kettenbruch diese Form haben:
\[
	1 - \cfrac{(nn-1)tt}{3 - \cfrac{(nn-4)tt}{5 - \cfrac{(nn-9)tt}{7 - \cfrac{(nn-16)tt}{9 - \mathrm{etc}}}}}
\]
aber die linke Seite wird nun in der Tat
\[
	\frac{nt\sqrt{-1}[(1+t\sqrt{-1})^n + (1-t\sqrt{-1})^n]}{(1+t\sqrt{-1})^n - (1-t\sqrt{-1})^n}
\]
sein, welche, weil die Imaginärteile nicht dagegen sprechen, gewiss einen reellen Wert haben muss, welchen wir also hier untersuchen finden wollen. Zu diesem Ziel wollen wir $t = \frac{\sin{\varphi}}{\cos{\varphi}}$ setzen, sodass $t = \tan{\varphi}$ ist; dann wird also
\[
	(1+t\sqrt{-1})^n = \frac{(\cos{\varphi} + \sqrt{-1}\sin{\varphi})^n}{(\cos{\varphi})^n} = \frac{\cos{n\varphi} + \sqrt{-1}\sin{n\varphi}}{(\cos{\varphi})^n}
\]
sein und auf ähnliche Weise
\[
	(1 - t\sqrt{-1})^n = \frac{(\cos{\varphi} - \sqrt{-1}\sin{\varphi})^n}{(\cos{\varphi})^n} = \frac{\cos{n\varphi} - \sqrt{-1}\sin{n\varphi}}{(\cos{\varphi})^n};
\]
Nachdem also diese Werte eingesetzt werden, wird die linke Seite
\[
	\frac{2n\sqrt{-1}\tan{\varphi}\cos{n\varphi}}{2\sqrt{-1}\sin{\varphi}} = \frac{n\tan{\varphi}\cos{n\varphi}}{\sin{n\varphi}} = \frac{n\tan{\varphi}}{\tan{n\varphi}}.
\]
\paragraph{§11}
Für $\tan{\varphi} = t$ gesetzt werden wir also den folgenden höchst bemerkenswerten Kettenbruch haben:
\[
	\cfrac{nt}{\tan{n\varphi}} = 1 - \cfrac{(nn-1)tt}{3 - \cfrac{(nn-4)tt}{5 - \cfrac{(nn-9)tt}{7 - \mathrm{etc}}}}
\]
welcher also auf folgende Weise ausgedrückt werden können wird:
\[
	\tan{n\varphi} = \cfrac{nt}{1 - \cfrac{(nn-1)tt}{3 - \cfrac{(nn-4)tt}{5 - \cfrac{(nn-9)tt}{7 - \mathrm{etc}}}}}
\]
welcher Ausdruck also sehr schön genutzt werden kann um die Tangenten vielfacher Winkel durch den Tangens des einfachen Winkels $t$ auszudrücken. Wenn so $n=2$ war, werden wir
\[
	\tan{2\varphi} = \frac{2t}{1-tt}
\]
haben. Wenn auf dieselbe Weise $n=3$ ist, wird
\[
	\tan{3\varphi} = \cfrac{3t}{1- \cfrac{8tt}{3-tt}} = \frac{3t - t^3}{1 - 3tt}
\]
sein. Hier offenbart sich ein höchst bemerkenswerter Fall, wann immer der Exponent $n$ unendlich klein angenommen wird; dann wird nämlich $\tan{n\varphi} = n\varphi$ sein; es wird also, indem man auf beiden Seiten durch $n$ teilt, diese Form entstehen
\[
	\varphi = \cfrac{t}{1 + \cfrac{tt}{3 + \cfrac{4tt}{5 + \cfrac{9tt}{7 + \mathrm{etc}}}}}
\] 
durch welchen Kettenbruch der Winkel selbst durch den Tangens $t$ ausgedrückt wird.
\paragraph{§12}
Wir wollen nun den Fall betrachten, in dem der Exponent $n$ unendlich groß angenommen wird, aber der Winkel $\varphi$ in der Tat unendlich klein und daher an auch sein Tangens $t$ unendlich klein ist, dennoch so, das $n\varphi = \theta$ ist und daher auch $nt=\theta$; dann werden wird also diesen Kettenbruch haben
\[
	\tan{\theta} = \cfrac{\theta}{1 - \cfrac{\theta\theta}{3 - \cfrac{\theta\theta}{5 - \cfrac{\theta\theta}{7 - \mathrm{etc}}}}}
\]
durch welche Formel aus einem gegebenen Winkel sein Tangens bestimmt werden können wird, welcher Ausdruck also als der reziproke des Vorhergehenden betrachtet werden kann.
\section*{Schlussfolgerung 2, \\ in dem ein verschwindender Winkel angenommen wird}
\paragraph{§13}
In diesem Fall wird der Kettenbruch
\[
	1 - \cfrac{zz}{3 - \cfrac{4zz}{5 - \cfrac{9zz}{7 - \cfrac{16zz}{9 - \mathrm{etc}}}}}
\]
sein. Für den linken Teil ist aber zu bemerken, dass
\[
	\frac{(1+z)^n - 1}{n} = \log{(1+z)}
\]
ist und daher
\[
	(1+z)^n = 1 + n\log{(1+z)};
\]
auf ähnliche Weise wird
\[
	(1-z)^n = 1 + n\log{(1-z)}
\]
sein, woher die linke Seite
\[
	\frac{nz[2 +n\log{(1+z)} + n\log{(1-z)}]}{n\log{(1+z)} - n\log{(1-z)}} = \frac{2z}{\log{\frac{1+z}{1-z}}}
\]
werden wird; daher werden wir also diese Form haben:
\[
	\cfrac{2z}{\log{\frac{1+z}{1-z}}} = 1 - \cfrac{zz}{3 - \cfrac{4zz}{5 - \cfrac{9zz}{7 - \cfrac{16zz}{9 - \mathrm{etc}}}}}
\]
und daher wird der Logarithmus selbst auf die folgende Weise ausgedrückt werden:
\[
	\log{\frac{1+z}{1-z}} = \cfrac{2z}{1 - \cfrac{zz}{3 - \cfrac{4zz}{5 - \mathrm{etc}}}}
\]
\section*{Schlussfolgerung 3, \\ in welcher der Exponent $n$ unendlich groß genommen wird}
\paragraph{§14}
Hier setze man also, damit der Kettenbruch einen endlichen Wert bekommt, was nur passieren kann, wenn die Größe $z$ unendlich klein gesetzt wird, $nz = v$, dass $z = \frac{v}{n}$ ist, und unser Kettenbruch wird
\[
	1 + \cfrac{vv}{3 + \cfrac{vv}{5 + \cfrac{vv}{7 + \cfrac{vv}{9 + \mathrm{etc}}}}}
\]
sein. Für die linke Seite ist aber bekannt, dass
\[
	\left( 1 + \frac{v}{n} \right)^n = e^v
\]
ist und auf die gleiche Art
\[
	\left( 1 - \frac{v}{n} \right) = e^{-v},
\]
die linke Seite wird also diese Form haben:
\[
	\frac{v(e^v + e^{-v})}{e^v - e^{-v}} = \frac{v(e^{2v} + 1)}{e^{2v} - 1};
\]
deswegen werden wir diesen bemerkenswerten Kettenbruch haben:
\[
	\frac{v(e^{2v} + 1)}{(e^{2v} - 1)} = 1 + \cfrac{vv}{3 + \cfrac{vv}{5 + \cfrac{vv}{7 + \cfrac{vv}{9 + \mathrm{etc}}}}}
\]
dessen transzendenter Wert auch auf diese Weise durch gewohnte Reihen beschafft werden kann:
\[
	\frac{1 + \frac{vv}{1\cdot 2} + \frac{v^4}{1\cdot 2\cdot 3\cdot 4} + \frac{v^6}{1\cdot 2\cdots 6} + \mathrm{etc}}{1 + \frac{vv}{1\cdot 2\cdot 3} + \frac{v^4}{1\cdot 2\cdots 5} + \frac{v^6}{1\cdot 2\cdots 7} + \mathrm{etc}}.
\]
\end{document}